\documentclass[a4paper,10pt]{amsart}

\usepackage[T1]{fontenc}
\usepackage[utf8]{inputenc}
\usepackage{amsmath,amssymb,amsthm,mathtools}
\usepackage{geometry}
\usepackage[colorlinks=true,allcolors=blue]{hyperref}

\theoremstyle{plain}
\newtheorem{theorem}{Theorem}[section]
\newtheorem{lemma}[theorem]{Lemma}
\newtheorem{proposition}[theorem]{Proposition}
\newtheorem{corollary}[theorem]{Corollary}
\newtheorem{problem}[theorem]{Problem}
\theoremstyle{definition}
\newtheorem{definition}[theorem]{Definition}

\theoremstyle{remark}
\newtheorem{remark}[theorem]{Remark}

\DeclareMathOperator{\ind}{ind}
\DeclareMathOperator{\Gal}{Gal}
\DeclareMathOperator{\Rad}{Rad}
\DeclareMathOperator{\End}{End}
\DeclareMathOperator{\Nrd}{Nrd}
\DeclareMathOperator{\Alt}{Alt}

\DeclareMathOperator{\id}{id}
\DeclareMathOperator{\car}{char}
\DeclareMathOperator{\GL}{GL}
\DeclareMathOperator{\Tr}{Tr}
\DeclareMathOperator{\res}{res}
\DeclareMathOperator{\ad}{ad}

\DeclareMathOperator{\rank}{rank}

\title[A Galois-type Springer theorem]
{A Galois-type Springer theorem for hermitian and generalized quadratic forms}

\author{A.-H. Nokhodkar}

\thanks{Department of Pure Mathematics, Faculty of Mathematical Sciences,
University of Kashan, Kashan, Iran.
Corresponding author.
\\
Email:
\href{mailto:a.nokhodkar@kashanu.ac.ir}
{a.nokhodkar@kashanu.ac.ir}.
\\
ORCID:
\href{https://orcid.org/0000-0003-0469-0072}
{0000-0003-0469-0072}.}


\begin{document}

\begin{abstract}
We prove that anisotropy of nonsingular hermitian and skew-hermitian forms
over central simple algebras with involution is preserved under scalar
extension to every finite Galois extension of odd degree coprime to the index
of the algebra. In characteristic two, the same conclusion holds for
generalized quadratic forms over algebras with involution of the first kind,
with no nonsingularity assumption on their polar forms. The main tool is an isometry-level descent theorem. For a cyclic extension of odd prime degree, every anisotropic nonsingular form over a division algebra with involution, of rank smaller than the extension degree and with Galois-invariant isometry class, descends to an anisotropic form over the base field. Unlike descent
results for Witt classes, this theorem descends the form itself. A
bounded-rank reduction, based on injectivity of the Witt-group restriction
map, derives the cyclic anisotropy theorem from this descent result; a
regularization argument  removes the nonsingularity assumption
on polar forms in characteristic two, and the Feit--Thompson theorem extends
the result to arbitrary finite Galois extensions of odd degree.
\end{abstract}

\maketitle
\section{Introduction}

A classical theorem of Artin--Springer states that an anisotropic quadratic
form remains anisotropic under every finite extension of odd degree
\cite{Springer1952,EKM}. It is natural to ask whether the same result holds
for hermitian forms over central simple algebras with involution. This
question is considerably more delicate. For quaternion division algebras with
an involution of the first kind, Parimala--Sridharan--Suresh proved an
odd-degree anisotropy theorem in characteristic different from $2$
\cite{PSS2001}. They also constructed an anisotropic unitary involution that
becomes isotropic over an odd-degree extension, and suggested a new
reformulation, in which the degree of the extension is coprime to
$2\ind(A)$. This reformulation was studied by Black and
Qu\'eguiner-Mathieu, who also established further cases for algebras of small
degree \cite{BlackQM2014}. As recalled there, results of Karpenko, Tignol and
Karpenko--Zhykhovich \cite{Karpenko2010,KZ2013} show that if anisotropy is
preserved under odd-degree extensions, then it is preserved under extension
to a generic splitting field of the algebra. Earlier, Lewis considered a
different extension map: for a field $K$ with a nontrivial involution and a
finite-dimensional $K$-algebra $A$ with an involution of the second kind
which is odd-dimensional modulo its radical, he proved that
$W(K)\to W(A)$ is injective, and he obtained a weaker Springer-type result
for central simple algebras of odd degree under a suitable restriction on the
involution \cite{Lewis1989}. In characteristic $2$, preservation of
anisotropy under odd-degree extensions was established for hermitian forms
over quaternion division algebras \cite{Nokhodkar2018}.

In characteristic $2$, a quadratic form is not determined in general by its polar form, so passing only to the associated bilinear or hermitian form may discard genuine quadratic information. This motivates working with generalized quadratic forms over algebras with involution, which retain this information and arise naturally in the theory of orthogonal groups and quadratic pairs \cite{et}.  In the broader theory, when the base field has characteristic different from $2$, or when the involution is unitary, the quadratic map is determined by its polar form; the genuinely additional case treated here is characteristic $2$ with an involution of the first kind.
For generalized quadratic forms over quaternion algebras, the
characteristic-free statement follows from the argument of \cite{PSS2001}
once one combines the excellence of function fields of conics
\cite{MT2016} with non-hyperbolic splitting of quadratic pairs
\cite{BD2015}; direct proofs have recently been given in
\cite{QMT2025,Nokhodkar2026}. The full form of question for arbitrary
finite odd-degree extensions remains open for general central simple
algebras with involution.

In this work we study the anisotropy of $(\pm1)$-hermitian and generalized quadratic forms under odd-degree Galois extensions. We treat all settings within one framework: nonsingular $(\pm1)$-hermitian forms for involutions of the first and second kinds in arbitrary characteristic, and generalized quadratic forms for involutions of the first kind in characteristic $2$. In the latter setting, we do not impose any nonsingularity assumption on the polar form.
Let $A$ be a central simple algebra over its center $K$, and let $\sigma$ be an involution on $A$. Put $F=K$ when $\sigma$ is of the first kind, and $F=K^\sigma$ when it is of the second kind; in the latter case $K/F$ is separable quadratic. Our main result states that, in each of the aforementioned settings, every anisotropic form over $(A,\sigma)$ remains anisotropic after scalar extension to $(A_L,\sigma_L)$ for every finite Galois extension $L/F$ of odd degree with $[L:F]$ coprime to $\ind(A)$ (see Corollary~\ref{csa-galois}). We note that when $\sigma$ is of the first kind, the coprimality condition is automatic, since $\ind(A)$ is a power of $2$; it is relevant only for $(\pm1)$-hermitian forms over central simple algebras with involution of the second kind.

All settings are handled by a unified proof scheme, rather than by independent arguments. The first step is to reduce the problem to the case where the underlying algebra is division, using Morita equivalence. In the quadratic setting, anisotropic generalized quadratic forms over a nonsplit division algebra, regardless of the degeneracy of their polar forms, are reduced to the nonsingular polar case by a regularization result (see Proposition~\ref{reg}), which is itself of independent interest. After this reduction, we consider cyclic extensions of odd prime degree. The key technical ingredient is a Galois descent statement for anisotropic forms of rank less than $p$ (Theorem~\ref{invariant}), which yields the cyclic case of prime degree (Theorem~\ref{cyclic}) and then, via the Feit--Thompson theorem \cite{FT}, the general Galois case (Theorem~\ref{galois}).

For Galois extensions of odd degree, descent of Witt classes is classical in
characteristic different from $2$: the Rosenberg--Ware theorem \cite{RW70}
identifies $W(F)$ with $W(L)^{\mathrm{Gal}(L/F)}$. Bayer-Fluckiger extended
this to hermitian forms over algebras with involution, and, by the same
method, extended a descent result of Rost \cite{Ros99} to arbitrary
odd-degree extensions \cite{BF09}. Injectivity of the restriction map
$W(A,\sigma)\to W(A_L,\sigma_L)$ in odd degree goes back to Bayer-Fluckiger
and Lenstra \cite{BL}; it holds for any finite-dimensional algebra with an
$F$-linear involution, in particular for involutions of the second kind.
Purely inseparable extensions in characteristic different from $2$ were
treated in \cite{BFM11}. However, descent of Witt classes alone does not imply descent of the hermitian forms themselves. As shown in \cite[Corollaries 4.2 and 4.5]{BF09}, such a descent follows from a Springer-type isotropy property for hermitian forms, namely, that isotropy over an odd-degree extension implies isotropy over the base field. This is precisely the open problem addressed in the present work. Thus, the descent result in [BF09] is conditional on the very property whose validity remains open, and cannot independently resolve this problem. In contrast, our descent result is established independently and provides a key ingredient in proving our main theorem. The regularization result in
Proposition~\ref{reg} then allows us to handle generalized quadratic forms
with arbitrary polar forms. Combining these results, our main theorem
establishes anisotropy preservation for finite Galois extensions of odd
degree.

\section{Preliminaries and conventions}\label{sec:notation}

Throughout this work, we fix a field $F$ of arbitrary characteristic and take $K$ to be either $F$
itself or a separable quadratic field extension of $F$. All algebras are
associative, unital and finite-dimensional.

A \emph{central simple $K$-algebra} is a $K$-algebra $A$ with centre exactly $K$
and no two-sided ideals other than $0$ and $A$. By Wedderburn's theorem \cite[Thm. III.4.1]{BO},
$A\simeq M_n(D)$ for a central division $K$-algebra $D$ that is unique
up to isomorphism; $D$ is said to be \emph{Brauer equivalent} to $A$. We put
$\deg A=\sqrt{\dim_KA}$ and $\ind A=\deg D$. The algebra $A$ is
\emph{split} if $A\simeq M_n(K)$, i.e.\ $\ind A=1$.

Let $A$ be a central simple $K$-algebra.
An \emph{involution} on $A$ is an additive map $\sigma:A\to A$ with
$\sigma(xy)=\sigma(y)\sigma(x)$ for all $x,y\in A$ and $\sigma^2=\id$.
Let $\sigma$ be an involution on $A$ satisfying
$K^\sigma=\{\lambda\in K:\sigma(\lambda)=\lambda\}=F$.
Since $\sigma$ preserves the centre $K$, the restriction $\sigma|_K$ is an
automorphism of $K$ of order $\le2$ with fixed field $F$. Hence exactly one of the
following holds:
\begin{itemize}
\item $K=F$ and $\sigma|_K=\id$: $\sigma$ is an involution \emph{of the
first kind}; it is $F$-linear.
\item $K/F$ is quadratic and $\sigma|_K$ is the non-trivial $F$-automorphism of
$K$: $\sigma$ is an involution \emph{of the second kind}; it is
$\sigma|_K$-semilinear.
\end{itemize}
In both cases $\sigma$ fixes $F$ and the pair $(A,\sigma)$ is called an \emph{algebra with involution}. If $\sigma$ is an involution of the first kind on $A$,
then $A\simeq A^{\mathrm{op}}$, where $A^{\mathrm{op}}$
denotes the opposite algebra. Consequently,
$2[A]=0$ in $\operatorname{Br}(F)$, so
$\exp(A)\mid 2$. By \cite[Cor.~2.8]{KMRT},
$\ind(A)$ is a power of $2$.

All modules considered in this work are finitely generated
projective right $A$-modules. The projectivity assumption
is automatic when $A$ is semisimple, in particular when
$A$ is central simple.
 For a field
extension $E/F$ such that $K_E=K\otimes_FE$ is a field, $A_E:=A\otimes_FE$ is a central simple $K_E$-algebra and
$\sigma_E:=\sigma\otimes\id$ is an involution of $A_E$ of the same kind, with
$K_E^{\sigma_E}=E$. For a module $V$ over $A$ we put $V_E=V\otimes_FE$, a module over $A_E$.
We use scalar extension only for such $E$.
\begin{lemma}\label{se}
Let $E/F$ be a finite extension of odd degree. Let $A$ be a central
simple $K$-algebra with an involution $\sigma$ satisfying $K^\sigma=F$,
and assume $\gcd([E:F],\ind A)=1$. Then:\\
(a) $K_E=K\otimes_FE$ is a field.\\
(b) $A_E$ is a central simple $K_E$-algebra with
$\ind A_E=\ind A$. In particular, if $A$ is a
division algebra, then so is $A_E$.\\
(c) $\sigma_E=\sigma\otimes\id$ is an involution of $A_E$ of the
same kind as $\sigma$, and $K_E^{\sigma_E}=E$.
\end{lemma}

\begin{proof}
(a) If $K=F$, then $K_E=E$. If $K/F$ is quadratic, then $K_E$ is a
field: a root of the irreducible quadratic polynomial of $K/F$ cannot
lie in $E$, because $[E:F]$ is odd.

(b) By \cite[Thm.~V.3.1(4)--(5)]{BO},
$\ind A_E\mid\ind A$ and
$\ind A\mid[E:F]\cdot\ind A_E$. Writing
$\ind A=k\cdot\ind A_E$, the second divisibility
gives $k\mid[E:F]$, hence $k\mid\gcd([E:F],\ind A)=1$,
so $k=1$ and $\ind A_E=\ind A$.

(c) Immediate from the definitions and from (a).
\end{proof}
Let $A$ be a central simple algebra over $K$ and let $V$ be a module over $A$.
A \emph{sesquilinear form} on $V$ is a biadditive map $f:V\times V\to A$ with
\[
 f(xa,yc)=\sigma(a)f(x,y)c\qquad(x,y\in V,\ a,c\in A).
\]
Put $f^*(x,y)=\sigma\bigl(f(y,x)\bigr)$, so that $f^{**}=f$. The form $f$ is
\emph{$\varepsilon$-hermitian} ($\varepsilon=\pm1$) if $f^*=\varepsilon f$, and
\emph{hermitian} if $\varepsilon=1$ (in characteristic $2$, the two signs coincide). The extension $f_E$ to $V_E$ is determined by
$f_E(x\otimes\lambda,y\otimes\mu)=f(x,y)\otimes\lambda\mu$, and
$(f_E)^*=(f^*)_E$.

For a sesquilinear form $b$ the \emph{adjoint map}
$\ad_b:V\to V^*=\operatorname{Hom}_A(V,A)$, defined by $\ad_b(x)=b(x,\cdot)$,
is additive and satisfies $\ad_b(xa)=\sigma(a)\ad_b(x)$
(here $V^*$ is a left $A$-module via $(a\varphi)(y)=a\varphi(y)$); in particular it is
$F$-linear. The form $b$ is \emph{nonsingular} if the adjoint map
$\ad_b$ is bijective.
The \emph{radical} of $b$ is defined by
$\Rad(b)=\ker(\ad_b)$.
 Since $V$ is finitely generated
projective, $(V_E)^*\simeq(V^*)_E$ and $\ad_{b_E}=\ad_b\otimes_FE$;
hence nonsingularity is preserved by scalar extension and
$\Rad(b_E)=\Rad(b)_E$.

Let $\car F=2$ and put
\[
 \Alt(A,\sigma)=\{a+\sigma(a):a\in A\},\quad
 \overline A=A/\Alt(A,\sigma),\quad \bar a=a+\Alt(A,\sigma).
\]
Then
$\sigma(c)\Alt(A,\sigma)\,c\subseteq\Alt(A,\sigma)$ for
$c\in A$, and
$\Alt(A_E,\sigma_E)=\Alt(A,\sigma)\otimes_FE$.

For a sesquilinear form $f$ on $V$ put
\[
 q_f(x)=\overline{f(x,x)}\in\overline A,\qquad b_f=f+f^* .
\]
The pair $q=(q_f,b_f)$ is called a \emph{generalized quadratic form} (with minimal form parameter
$\Alt(A,\sigma)$) \emph{represented by} $f$, and we write $b_q=b_f$,
see \cite[Ch.~I]{knus}. A generalized quadratic form thus consists of both $q$ and its polar form $b_q$, and
no nonsingularity of $b_q$ is assumed. One checks directly that for $x,y\in V$,
$a,c\in A$, $\lambda\in F$:
\[
 q(x+y)=q(x)+q(y)+\overline{b_q(x,y)},\quad
 q(xa)=\overline{\sigma(a)f(x,x)a},\quad
 b_q(xa,yc)=\sigma(a)b_q(x,y)c,
\]
and $b_q^*=b_q$, $q(x\lambda)=\lambda^2q(x)$. The scalar extension is
$q_E=(q_{f_E},b_{f_E})$; it is independent of the choice of $f$ by
Lemma~\ref{quadratic}(a) below.

An $\varepsilon$-hermitian form $h$ (resp.\ a generalized quadratic form $q$) is called
\emph{isotropic} if $h(x,x)=0$ (resp.\ $q(x)=0$) for some $x\ne0$, and
\emph{anisotropic} otherwise. A submodule $U\subseteq V$ is \emph{totally isotropic} if $h|_{U\times U}=0$
(resp.\ $q|_U=0$ and $b_q|_{U\times U}=0$, i.e.\ the generalized quadratic form
vanishes on $U$);
compare \cite[(3.7)]{knus}, where $U\subseteq U^\perp$ for the associated
hermitian form.

\begin{lemma}\label{quadratic}
(a) Let $\car F=2$. Two sesquilinear forms represent the same generalized quadratic form if and only if there exists a sesquilinear form $r$ on $V$ such that
$f-g=r+r^*$.
 In particular $(q_f,b_f)_E:=(q_{f_E},b_{f_E})$ is well defined, and $q|_U=0=b_q|_U$ if and only if $f|_U\in\{r+r^*: r \text{ sesquilinear on } U\}$.
\\
(b) For every $\varepsilon$-hermitian form
($\varepsilon=\pm1$) in any characteristic,
and for generalized quadratic forms in
characteristic $2$, a form is isotropic
if and only if it has a nonzero totally
isotropic submodule.
\end{lemma}

\begin{proof}
(a) If $f-f'=r+r^*$, then $b_f=b_{f'}$ and
\[
f(x,x)-f'(x,x)=r(x,x)+\sigma(r(x,x)).
\]

Conversely, let $t=f-f'$. Then $t=t^*$ and
$t(x,x)\in\Alt(A,\sigma)$ for every $x\in V$.

First suppose that $V=A^n$, and let
$e_1,\ldots,e_n$ denote the standard basis of
the free right $A$-module $A^n$.
Set
$T_{ij}=t(e_i,e_j)$, $1\leq i,j\leq n$.
Since $t=t^*$, we have
$T_{ij}=\sigma(T_{ji})$.
Moreover, since $T_{ii}=t(e_i,e_i)\in
\Alt(A,\sigma)$, for each $i$ we can choose
$c_i\in A$ such that
$T_{ii}=c_i+\sigma(c_i).$

Define a sesquilinear form $r$ on $A^n$ by
specifying its matrix $R=(R_{ij})$ as follows:
\[
R_{ij}=
\begin{cases}
T_{ij},&i<j,\\
c_i,&i=j,\\
0,&i>j.
\end{cases}
\]
Then, by construction, $t=r+r^*$.

For a general finitely generated projective module
$V$, choose a complement $W$ such that
$
V\oplus W\simeq A^n
$
for some $n$. Let
$\pi:V\oplus W\longrightarrow V$
be the canonical projection onto $V$ along $W$.
The form
$
\widetilde t(x,y)=t(\pi(x),\pi(y))$,
$ x,y\in V\oplus W,
$
is sesquilinear and hermitian, and
$
\widetilde t(x,x)\in\Alt(A,\sigma).
$
Identifying $V\oplus W$ with $A^n$, the free-module
case gives a sesquilinear form $\widetilde r$
such that
$
\widetilde t=\widetilde r+\widetilde r^*.
$
Restricting $\widetilde r$ to $V$ yields the desired
sesquilinear form $r$ on $V$ satisfying
$
t=r+r^*.
$
The last assertion follows by applying the first
part to $f|_U$ and the zero form.

(b) If $h(x,x)=0$ then $xA$ is totally isotropic. If $q(x)=0$ then $f(x,x)\in\Alt$, so
$b_q(x,x)=0$ in characteristic~$2$, and $xA$ is totally isotropic. The converse is clear.
\end{proof}

In characteristic $2$, generalized quadratic forms with
involutions of the second kind are completely determined
by their polar hermitian forms. Indeed, in this case
$
\operatorname{Sym}(A,\sigma)
=
\Alt(A,\sigma),
$
see \cite[Prop.~2.17]{KMRT}. If $f$ and $g$ have
the same polar form, then $t=f-g$ satisfies $t=t^*$,
and hence $t(x,x)\in\operatorname{Sym}(A,\sigma)
=\Alt(A,\sigma)$ for every $x\in V$.
Thus $q_f=q_g$. Consequently, the genuinely quadratic
case not determined by polarization occurs only for
involutions of the first kind. We therefore restrict
the generalized quadratic form case below to $K=F$.

\begin{problem}\label{csa}
Let $E/F$ be finite extension of odd degree, and let $A$ be a central simple algebra
over $K$, with
$\gcd([E:F],\ind A)=1$. Let $\sigma$ be an involution of $A$ with
$K^\sigma=F$ and let $V$ be an $A$-module. Assume one of the following:
\\
(H) $\varphi=h$ is a nonsingular
$\varepsilon$-hermitian form over $(A,\sigma)$,
where $\varepsilon=\pm1$, with no restriction on $\car F$.
\\
(Q) $\car F=2$, $K=F$, and $\varphi=q$ is a generalized
quadratic form with minimal form parameter
$\Alt(A,\sigma)$ (no condition on $b_q$).

Does the anisotropy of $\varphi$ over $F$ imply the anisotropy of its scalar extension $\varphi_E$ over $E$?
\end{problem}

Note that Problem~\ref{csa} in setting (H) is a hermitian reformulation of question~(i$'$) of \cite[\S0]{BlackQM2014}.
Also, if $K=F$, then $\exp A\le2$, so $\ind A$ is a power of $2$ and the index condition is automatic.

\section{Morita equivalence and reduction of problem}
The aim of this section is to reduce Problem \ref{csa} to the case where $A$ is a division algebra in both settings, $\varepsilon=1$ (i.e., only hermitian forms need to be considered) in setting (H), and the generalized quadratic form $q$ is nonsingular in setting (Q). We continue with the notation of the preceding section.

\begin{lemma}\label{md}
Let $D$ be the central division $K$-algebra Brauer equivalent to $A$. There are
an involution $\tau$ of $D$, a right $D$-module $M\ne0$ and a nonsingular
$\varepsilon_0$-hermitian form $b_0$ on $M$ ($\varepsilon_0=\pm1$) such that
$A\simeq\operatorname{End}_D(M)$ and $\sigma=\sigma_{b_0}$ is the adjoint involution of $\operatorname{End}_D(M)$ with respect to $b_0$.
Moreover $\tau|_K=\sigma|_K$, $K^\tau=F$ and $\ind D=\ind A$.
These identifications are compatible
with scalar extension to every field extension
$E/F$ for which $K_E$ is a field.
\end{lemma}

\begin{proof}
The first assertion is the standard description of
involutions as adjoint involutions
\cite[Thm.~4.2]{KMRT}.
For $\lambda\in K$ acting on $M$ we have
$b_0(m\lambda,m')=\tau(\lambda)b_0(m,m')=b_0(m,m'\tau(\lambda))$, so the adjoint of
$\lambda$ is $\tau(\lambda)$, i.e.\ $\sigma|_K=\tau|_K$. Since $D$ is Brauer equivalent to $A$, their indices
coincide. Moreover, $\tau|_K=\sigma|_K$, so
$K^\tau=F$. The compatibility with scalar extension
follows from
$A_E\simeq\operatorname{End}_{D_E}(M_E)$ and
$\sigma_E=\sigma_{(b_0)_E}$.
\end{proof}

\begin{proposition}\label{transfer}
With the notation of Lemma~\ref{md}, let $W=V\otimes_AM$ and, for a
sesquilinear form $k$ on $V$, define
\[
 \widehat k(v\otimes m,v'\otimes m')=b_0\bigl(m,k(v,v')m'\bigr).
\]
Then $k\mapsto\widehat k$ is a bijection $\operatorname{Sesq}_A(V)\to\operatorname{Sesq}_D(W)$
with $(\widehat k)^*=\varepsilon_0\widehat{k^*}$, compatible with restriction to
submodules and with scalar extension to $E$. Consequently:
in (H), $h$ goes to a nonsingular
$(\varepsilon_0\varepsilon)$-hermitian form
$h_D=\widehat h$ over $(D,\tau)$; in (Q),
$q$ goes to a generalized quadratic form $q_D$
over $(D,\tau)$, with no nonsingularity assumption
on its polar form. In both cases anisotropy
over $F$ and over $E$ is preserved, i.e.,
$\varphi$ and $\varphi_E$ are anisotropic if and
only if $\varphi_D$ and $(\varphi_D)_E$ are.
\end{proposition}

\begin{proof}
This is the usual Morita theory of sesquilinear forms \cite[Ch.~I, \S9]{knus};
we sketch it. Since $V$ is a finitely generated projective
$A$-module, it is a direct summand of $A^r$
for some $r\geq 1$. By additivity and passage
to direct summands, it therefore suffices to
consider the case $V=A$, where
$\operatorname{Sesq}_A(A)\simeq A$ and $\gamma\mapsto b_0(\cdot,\gamma\,\cdot)$ is a bijection
$A\to\operatorname{Sesq}_D(M)$, because $b_0$ is nonsingular. Using
$\tau(b_0(x,y))=\varepsilon_0b_0(y,x)$ and $b_0(\beta m,m')=b_0(m,\sigma(\beta)m')$, we get
\[
 (\widehat k)^*(v\otimes m,v'\otimes m')
 =\varepsilon_0\,b_0\bigl(k(v',v)m,m'\bigr)=\varepsilon_0\,\widehat{k^*}(v\otimes m,v'\otimes m').
\]
Compatibility with restriction and with extension to $E$ (use
$A_E\simeq\operatorname{End}_{D_E}(M_E)$) is clear from the formula; nonsingularity
of $\widehat h$ in (H) is \cite[Ch. I, Thm. 9.3.5]{knus}.

Morita equivalence gives a bijection $U\mapsto U\otimes_AM$ between submodules of $V$ and of $W$,
preserving nonzero modules. In (H), $\widehat h|_{U\otimes M}=0$ if and only if $h|_U=0$, since
$b_0(m,\beta m')=0$ for all $m,m'$ forces $\beta=0$ ($b_0$ nonsingular, $A$ acts faithfully on $M$).
In (Q), $\varepsilon_0=1$, so transport commutes with $r\mapsto r+r^*$ and maps
$\{r+r^*\}$ bijectively onto the corresponding set for $U\otimes M$; by
Lemma~\ref{quadratic}(a), $U$ is totally isotropic if and only if $U\otimes M$ is.
Now apply Lemma~\ref{quadratic}(b); the same argument works over $E$.
\end{proof}
\begin{lemma}
\label{tsi}
Let $E/F$ be a finite extension of odd degree in
characteristic two.
Suppose that $\tau$ is of the first kind.
Let $q$ be a generalized quadratic form over $(D,\tau)$,
and let $R\subseteq V$ be a right $D$-subspace satisfying
$b_q|_{R\times R}=0$.
If $q|_R$ is anisotropic, then $(q|_R)_E$
is anisotropic.
\end{lemma}

\begin{proof}
The restriction $q|_R$ is totally singular.
Moreover, $E/F$ is separable and $D_E$ is a division
algebra by Lemma \ref{se}. The assertion therefore follows from
\cite[Cor.~4.2]{TS}.
\end{proof}

\begin{corollary}
\label{outside}
Assume that $\car F=2$.
Let $q$ be anisotropic and let $E/F$ be a finite
extension of odd degree.
Every nonzero isotropic vector of $q_E$ lies outside
$\Rad((b_q)_E)$.
\end{corollary}

\begin{proof}
Apply Lemma~\ref{tsi} to
$R=\Rad(b_q)$.
Since $q$ is anisotropic, its restriction to $R$
is anisotropic. Moreover,
$\Rad((b_q)_E)=\Rad(b_q)_E$.
The assertion follows.
\end{proof}

We now settle singular quadratic forms.
\begin{proposition}
\label{reg}
Assume that $\car F=2$.
Let $(D,\tau)$ be a central division $F$-algebra
with involution of the first kind, and suppose
that $D\ne F$.
Let $q$ be an anisotropic generalized quadratic form
on a right $D$-vector space $V$ of dimension $n>0$,
with arbitrary polar form.
Put $F_t=F((t))$.
There exists a generalized quadratic form $q_t$
on $V_{F_t}$ such that:
\begin{enumerate}
\item $q_t$ has nonsingular polar form and is
anisotropic over $F_t$;
\item for every finite extension $E/F$ of odd degree,
if $q_E$ is isotropic, then
$(q_t)_{E((t))}$ is isotropic.
\end{enumerate}
\end{proposition}

\begin{proof}
Fix a right $D$-basis of $V$.
Since $D\ne F$, the involution $\tau$ is nontrivial:
otherwise, as $\tau$ is an anti-automorphism, $D$
would be commutative and hence equal to $F$.
Choose $a\in D$ such that
$c=a+\tau(a)\ne0$.
Then $\tau(c)=c$, and $c$ is invertible.
Define the generalized quadratic form $\psi$ by
\[
\psi(x_1,\ldots,x_n)
=
\overline{\sum_{i=1}^n\tau(x_i)a x_i}.
\]
A sesquilinear form representing $\psi$ is given by
$
g(x,y)=\sum_{i=1}^n\tau(x_i)a y_i.
$
Its polar form is
$
b_\psi(x,y)=\sum_{i=1}^n\tau(x_i)c y_i.
$
Thus its polar Gram matrix is $C=cI_n$, which is
invertible.
Let $B$ be the polar Gram matrix of $q$, and set
$q_t=q_{F_t}+t\psi_{F_t}$.
The polar Gram matrix of $q_t$ is $B+tC$.
Its reduced norm is a polynomial in $t$ whose leading
coefficient is $\Nrd(C)\ne0$.
Consequently, $B+tC$ is invertible over $D((t))$,
and the polar form of $q_t$ is nonsingular.
Coefficientwise, we have
\[
D_{F_t}\simeq D((t)),
\qquad
D_{F_t}/\Alt(D_{F_t},\tau_{F_t})
\simeq \overline D((t)),
\]
where
$\overline D=D/\Alt(D,\tau)$.
Since $D$ is a division algebra, so is $D((t))$.

We next show that $q_t$ is anisotropic.
For any nonzero $x\in V_{F_t}$, write
\[
x=t^m u,\qquad
u=x_0+tx_1+t^2x_2+\cdots,\quad x_0\ne0.
\]
Then
$q_t(x)=t^{2m}q_t(u)$,
and the constant term of $q_t(u)$ is $q(x_0)$,
which is nonzero by the anisotropy of $q$.
Hence $q_t$ is anisotropic over $F_t$.

It remains to prove the second assertion.
Note that $D_E$ is a division algebra by
Lemma~\ref{se}(b): since $\car F=2$
and $\tau$ is of the first kind, $\ind D$ is a power of
$2$, while $[E:F]$ is odd, so the two are coprime.
Choose $0\ne v\in V_E$ such that $q_E(v)=0$.
By Corollary~\ref{outside}, $v$ lies outside the radical of the
polar form. Since
$D_E$ is division, there exists $w\in V_E$ such that $b_E(v,w)=1$.
Let $f$ be a sesquilinear form representing $q$,
and let $g$ be the sesquilinear form representing
$\psi$ as above. Then $q_t$ is represented by
$f_{F_t}+tg_{F_t}$, and its polar form is
$b_t=b_q+t b_\psi$.

We construct a sequence $v_j\in V_E[[t]]$ satisfying
\[
v_j\equiv v\pmod t,\qquad
(q_t)_{E((t))}(v_j)\equiv0\pmod{t^j}.
\]
Start with $v_1=v$.
Suppose that $v_j$ has been constructed.
Let $\overline{d_j}\in\overline D_E$ be the
coefficient of $t^j$ in $(q_t)_{E((t))}(v_j)$,
and choose $d_j\in D_E$ representing it.
Thus
\[
(q_t)_{E((t))}(v_j)
\equiv t^j\overline{d_j}\pmod{t^{j+1}}.
\]
Put
$v_{j+1}=v_j+t^j w d_j$.
Since $v_j\equiv v\pmod t$ and
$b_t\bmod t=b_E$, we have
\[
b_t(v_j,w d_j)
\equiv b_E(v,w)d_j
=d_j\pmod t.
\]
Using the quadratic identity and $2j\ge j+1$, we obtain
\begin{align*}
(q_t)_{E((t))}(v_{j+1})
&=(q_t)_{E((t))}(v_j)
+t^j\overline{b_t(v_j,w d_j)}
+t^{2j}(q_t)_{E((t))}(w d_j)\\
&\equiv t^j\overline{d_j}
+t^j\overline{d_j}
\equiv0\pmod{t^{j+1}},
\end{align*}
where the last equality uses $\car F=2$.

The sequence $(v_j)$ converges $t$-adically to an element
$\widetilde v\in V_E[[t]]\subseteq V_{E((t))}$.
We have $\widetilde v\equiv v\ne0\pmod t$.
For every $j$, write
$\widetilde v-v_j=t^jz_j$, where $z_j\in V_E[[t]]$.
By the quadratic identity,
\[
(q_t)_{E((t))}(\widetilde v)
=(q_t)_{E((t))}(v_j)
+t^{2j}(q_t)_{E((t))}(z_j)
+t^j\overline{b_t(v_j,z_j)}.
\]
Consequently,
\[
(q_t)_{E((t))}(\widetilde v)
\equiv(q_t)_{E((t))}(v_j)
\equiv0\pmod{t^j}.
\]
Thus $(q_t)_{E((t))}(\widetilde v)=0$, and
$(q_t)_{E((t))}$ is isotropic.
\end{proof}

In each setting the \emph{division case} is the following assertion, required
over every base field $F$ of the relevant characteristic:
\begin{itemize}
\item[(H)] Problem~\ref{csa} has an affirmative answer when $A$ is a division
algebra and $\varepsilon=+1$, i.e.\ for nonsingular \emph{hermitian} forms over
a division algebra with involution of the first or second kind;
\item[(Q)] Problem~\ref{csa} has an affirmative answer when $A\ne F$ is a
division algebra and $b_q$ is nonsingular.
\end{itemize}

\begin{theorem}\label{mor}
If the division case holds in setting (H) (resp.\ (Q)), then Problem~\ref{csa}
has an affirmative answer in that setting for all $A$.
\end{theorem}

\begin{proof}
Let $\varphi$ be anisotropic, and take $D,\tau,M,b_0$ as in
Lemma~\ref{md}. By Proposition~\ref{transfer}, $\varphi_D$ is
anisotropic on $W=V\otimes_AM$, and $\varphi_E$ is anisotropic if and only if $(\varphi_D)_E$ is.
Also $K^\tau=F$ and $\gcd([E:F],\ind D)=\gcd([E:F],\ind A)=1$.
Since $D$ is a division algebra, Lemma~\ref{se}(b)
shows that $D_E$ is a division algebra. We may
assume $W\neq0$.

\smallskip\noindent\emph{Setting (H).}
Put $\varepsilon_1=\varepsilon_0\varepsilon$, so $h_D$ is a nonsingular
$\varepsilon_1$-hermitian form over $(D,\tau)$.
If $\varepsilon_1=1$ or $\car F=2$ (where $\varepsilon_1=1$ by
convention), the division case applies directly to $h_D$.
Let $\car F\ne2$ and $\varepsilon_1=-1$, so
$\tau(h_D(x,y))=-h_D(y,x)$.
If $\tau=\id$, then $D$ is commutative, so $D=K=F$, and
$h_D(x,x)=-h_D(x,x)$ gives $h_D(x,x)=0$; thus $h_D$ is isotropic on
$W\ne0$, a contradiction.
So $\tau\ne\id$; pick $c\in D$ with $\tau(c)\ne c$ and put
$a=c-\tau(c)\in D^\times$, so $\tau(a)=-a$.
Let $\tau'=\operatorname{Int}(a)\circ\tau$, i.e.\ $\tau'(x)=a\tau(x)a^{-1}$.
Then $\tau'$ is an anti-automorphism, and
$\tau'^2=\operatorname{Int}(a)\operatorname{Int}(\tau(a)^{-1})
=\operatorname{Int}(a)\operatorname{Int}(a^{-1})=\id$.
Hence $\tau'$ is an involution of $D$ with $\tau'|_K=\tau|_K$, so $K^{\tau'}=F$
and $\tau'$ has the same kind as $\tau$.
Put $h'=a\,h_D$.
Then
$h'(x\alpha,y\beta)=a\tau(\alpha)h_D(x,y)\beta=\tau'(\alpha)h'(x,y)\beta$,
$h'$ is nonsingular, and
\[
\tau'\bigl(h'(x,y)\bigr)=a\,\tau\bigl(h_D(x,y)\bigr)\tau(a)a^{-1}
=-a\,\tau\bigl(h_D(x,y)\bigr)=a\,h_D(y,x)=h'(y,x),
\]
so $h'$ is a nonsingular hermitian form over $(D,\tau')$.
Since $h'(x,x)=0$ if and only if $h_D(x,x)=0$ (over $F$ and over $E$), $h'$ and $h'_E$ are
anisotropic if and only if $h_D$ and $(h_D)_E$ are.
Moreover $K_E^{\tau'_E}=E$ and $\ind D$ is unchanged, so the
division case applies to $h'$ and shows that $h'_E$, hence $\varphi_E$, is
anisotropic.

\smallskip\noindent\emph{Setting (Q).} Here $K=F$. If $D=F$, then $\Alt(F,\id)=0$
and $q_D$ is an ordinary quadratic form, possibly singular; $(q_D)_E$ is anisotropic
by Springer's theorem, which holds for arbitrary quadratic forms in characteristic $2$
\cite[Cor.~18.5]{EKM}. Let $D\ne F$ and suppose $(q_D)_E$ is isotropic. Identify
$W=D^m$ and let $Q$ be as in Proposition~\ref{reg}, a form
over $(D_L,\tau_L)$, $L=F((t))$. Now $D_L=D((t))$ is a central division
$L$-algebra with involution of the
first kind, $D_L\ne L$, and $[E((t)):F((t))]=[E:F]$ is odd. Since $Q$ is anisotropic with
nonsingular polar form, the division case gives that $Q_{E((t))}$ is anisotropic, contradicting
Proposition~\ref{reg}. Hence $(q_D)_E$, and thus $q_E$, is anisotropic.
\end{proof}
\begin{remark}\label{new}
The proof of Theorem~\ref{mor} works verbatim if Problem~\ref{csa} and the division case are restricted to finite Galois extensions $E/F$ of odd degree (over every base
field of the relevant characteristic): the only extensions that occur in the proof are $E/F$ and $E((t))/F((t))$, and the latter is again Galois of degree $[E:F]$.
\end{remark}

\section{New setting for forms over division algebras}
\label{sec:forms}

We continue to assume that $F$ is a field of arbitrary characteristic and take $K$
to be either $F$ itself or a separable quadratic extension
of $F$. We also fix $D$ as a central division $K$-algebra equipped
with an involution $\tau$ satisfying $K^\tau=F$.

We write
$
\deg D=2^s m,
$
where $s\geq 0$ and $m$ is odd.
If $K=F$, then $\tau$ is an involution of the first
kind. In this case, the existence of $\tau$ implies
that $[D]$ has exponent dividing $2$. Since $D$ is
a division algebra, its index equals its degree and
is a power of $2$. Consequently, $m=1$.
If $K/F$ is quadratic, then $\tau$ is a unitary
involution, and no further restriction on the odd
factor $m$ is imposed.
Let $V$ be a finite-dimensional right $D$-vector space.
\begin{lemma}\label{polar}
Let $\car F=2$, $K=F$ and let $q$ be a generalized quadratic form over $(D,\tau)$.
Then the polar form $b_q$ is uniquely determined by the
quadratic map $q:V\to D/\Alt(D,\tau)$.
If $g$ is a sesquilinear form representing $q$, then
\[
q(x)=\overline{g(x,x)},\qquad
b_q(x,y)=g(x,y)+\tau(g(y,x)).
\]
In particular,
$b_q(x,x)\in\Alt(D,\tau)$, and
$q(x)=0$ implies $b_q(x,x)=0$.
\end{lemma}

\begin{proof}
See \cite[Prop. 1.1]{et} and its proof.
\end{proof}

Henceforth we write $q$ for the pair $(q,b_q)$.

In view of Theorem~\ref{mor}, it suffices to consider the
following two settings (the case $D=F$ in the quadratic setting was settled
in the proof of that theorem, by Springer's theorem).

\smallskip
\noindent\textbf{Hermitian setting.}
$\varphi=h$ is a nonsingular hermitian form over $(D,\tau)$.

\smallskip
\noindent\textbf{Quadratic setting.}
$\car F=2$, $K=F$, $D\ne F$, and
$\varphi=q$ is a generalized quadratic form over $(D,\tau)$
whose polar form $b_q$ is nonsingular.

\smallskip
\noindent \textbf{Common notation.}
We use $\varphi$ for a form of either type.
Its \emph{associated hermitian form} is $b_\varphi=h$ if $\varphi=h$ is hermitian and $b_\varphi=b_q$ if $\varphi=q$ is quadratic.
An isometry of $\varphi$ always means an isometry of
the form itself, not merely of $b_\varphi$.
In both settings $b_\varphi$ is nonsingular
by assumption.

\begin{definition}
An \emph{admissible odd-degree extension} of $F$ is a
finite extension $E/F$ such that
$[E:F]$  is odd  and $\gcd([E:F],m)=1$,
where $\deg D=2^s m$ with $s\geq 0$ and $m$ odd.
\end{definition}

As already observed, if $\tau$ is of the first kind, then $m=1$; hence every
finite odd-degree extension of $F$ is admissible.
For every admissible extension $E/F$, the algebra $K_E$ is a field and
$D_E$ remains a division algebra by
Lemma~\ref{se}(a)--(b): admissibility gives
$[E:F]$ odd and $\gcd([E:F],m)=1$, hence
$\gcd([E:F],\ind D)=\gcd([E:F],2^sm)=1$.

\medskip
\noindent \textbf{Witt decompositions}
For a generalized quadratic form over $(D,\tau)$,
a $D$-subspace $U$ is totally isotropic if and only if
$q|_U=0$.
Indeed, by Lemma~\ref{polar}, the condition $q|_U=0$
also implies $b_q|_{U\times U}=0$:
for all $x,y\in U$ and $a\in D$, polarization gives
\[
0=q(xa+y)-q(xa)-q(y)
=\overline{\tau(a)b_q(x,y)}.
\]
Thus $\tau(a)b_q(x,y)\in\Alt(D,\tau)$ for all $a\in D$.
If $b_q(x,y)\ne0$, this implies
$\Alt(D,\tau)=D$, a contradiction since
$\Alt(D,\tau)=\operatorname{im}(1+\tau)$
and $(1+\tau)^2=0$ in characteristic $2$.
Hence $b_q(x,y)=0$.
Thus the latter condition in the definition of total
isotropy in Section~\ref{sec:notation} is automatic here.
We write $i(\varphi)$ for the maximum dimension of
such a subspace.
The form $\varphi$ is called \textit{metabolic} if it is nonsingular and $i(\varphi)=\frac{1}{2}\rank \varphi$, where
$\rank \varphi=\dim_D V$.
A totally isotropic submodule $\Lambda\subseteq V$ satisfying
$\dim_D \Lambda=\frac{1}{2}\dim_D V$
is called a \textit{Lagrangian}.
Thus, a nonsingular form is metabolic if and only if
it admits a Lagrangian.
In the quadratic setting, $\varphi$ is metabolic if and only if it is hyperbolic (see \cite[Ch. I, Rem. 5.5.4]{knus}).

\begin{proposition}
\label{witt}
Every nonsingular form $\varphi$ under consideration
has a decomposition
$\varphi\simeq \varphi_{\mathrm{met}}\perp\varphi_{\mathrm{an}}$,
where $\varphi_{\mathrm{an}}$ is anisotropic and $\varphi_{\mathrm{met}}$ is metabolic with $\rank \varphi_{\mathrm{met}}=2i(\varphi)$.
The anisotropic representative of a Witt class is
unique up to isometry.
\end{proposition}

\begin{proof}
For hermitian forms, the decomposition and uniqueness follow from
\cite[Ch. I, Props. 6.1.1 and 6.1.4]{knus}.
For generalized quadratic forms, the result follows from
\cite[Ch. I, Thm. 6.5.1]{knus}.
\end{proof}

Write $W_{\mathrm H}(D,\tau)$ for the Witt group of
nonsingular hermitian forms modulo metabolic forms.
In the quadratic setting, write $W_{\mathrm Q}(D,\tau)$
for the Witt group of nonsingular generalized quadratic
forms modulo hyperbolic forms.
When the type of form is indicated by $\varphi$, we use
$W_\varphi(D,\tau)$ as shorthand for $W_{\mathrm H}(D,\tau)$
in the hermitian setting and for $W_{\mathrm Q}(D,\tau)$
in the quadratic setting. Thus the subscript $\varphi$
indicates only the type of form, not a particular form
or its isometry class. We write
\[
\res_{E/F}:W_\varphi(D,\tau)\longrightarrow
W_\varphi(D_E,\tau_E)
\]
for the scalar-extension homomorphism induced by
$\varphi\mapsto\varphi_E$.
\begin{proposition}
\label{injectivity}
For every admissible odd-degree extension $E/F$,
the map
$\res_{E/F}$
is injective.
\end{proposition}

\begin{proof}
In the hermitian setting, the assertion is \cite[Prop. 1.2]{BL},
viewing $D$ as an $F$-algebra: the involution $\tau$ is
$F$-linear, including when it is of the second kind.
In the quadratic setting, let $\psi$ be an arbitrary
nonsingular generalized quadratic form over $(D,\tau)$,
and denote its scalar extension to $D_E$ by $\psi_E$.
Put $t=\Tr_{E/F}$ and
$\lambda=\id_D\otimes t:D_E\to D$, where $\Tr_{E/F}$ is the trace of $E/F$.
The map $\lambda$ is $D$-bilinear and commutes with the
involutions, so
$\lambda(a+\tau_E(a))=\lambda(a)+\tau(\lambda(a))$.
Thus it maps $\Alt(D_E,\tau_E)$ into
$\Alt(D,\tau)$. Restricting scalars from
$D_E$ to $D$ and applying $\lambda$ to the values of
the form therefore defines a well-defined transfer
$t_*$ on generalized quadratic forms.
This transfer preserves nonsingularity: for a
nonsingular form $Q$ over $D_E$ and $x\ne0$,
one can choose $y$ with $b_Q(x,y)=1$, since
$D_E$ is division. Hence
$\lambda(b_Q(x,y))=t(1)=1$.
It also preserves hyperbolicity: if $\Lambda$ is a
Lagrangian of $Q$, then, viewed as a $D$-subspace
by restriction of scalars, it remains totally
isotropic for $t_*Q$ and has half the $D$-rank
of the underlying space. Since $t_*Q$ is
nonsingular, it is metabolic and hence hyperbolic.
Therefore $t_*$
induces a homomorphism on Witt groups.

Consider the trace form
$\beta(a,b)=t(ab)$ on $E$.
Since $\car F=2$ and $[E:F]$ is odd,
the extension $E/F$ is separable, so the trace pairing
$\beta(a,b)=t(ab)$ is nonsingular.
Moreover, $t(1)=1$, and hence
$E=F\cdot1\perp\ker t$.
For every $a\in\ker t$,
\[
\beta(a,a)=t(a^2)=t(a)^2=0.
\]
Thus $\beta|_{\ker t}$ is nonsingular alternating.
Put $n_E=([E:F]-1)/2$.
Choose a symplectic basis
$e_1,f_1,\ldots,e_{n_E},f_{n_E}$ of $\ker t$, so
$\beta(e_i,e_j)=\beta(f_i,f_j)=0$ and $\beta(e_i,f_j)=\beta(f_j,e_i)=\delta_{ij}$.
On the summand corresponding
to a hyperbolic plane with basis $e_i,f_i$, the transferred
quadratic map on $V\oplus V$ is
$(x,y)\mapsto\overline{b_\psi(x,y)}$, and its polar
form has block matrix
$\left(\begin{smallmatrix}0&b_\psi\\b_\psi&0\end{smallmatrix}\right)$.
Thus $V\oplus0$ is a Lagrangian, so this summand is
hyperbolic. The $F\cdot1$ summand contributes $\psi$.
Consequently, for the form $\psi$ fixed above,
$t_*(\psi_E)\simeq\psi\perp\psi_h$,
where $\psi_h$ is hyperbolic. Hence,
$
t_*\circ\res_{E/F}=\id
$
on $W_{\mathrm Q}(D,\tau)$, proving injectivity.
\end{proof}

\section{Index, embeddings, and centralizer norms}
\label{sec:algebra}
Let $B$ be a central simple $k$-algebra and
$E/k$ a finite field extension.
We use the standard index divisibilities
\[
    \ind(B_E)\mid\ind(B),
    \quad
    \ind(B)\mid[E:k]\cdot\ind(B_E);
\]
see \cite[Thm. V.3.1(4)--(5)]{BO}.
If, in addition, $E$ is a subfield of $B$,
the centralizer theorem \cite[Cor. III.5.3]{BO} gives
$B_E\simeq M_{[E:k]}(C_B(E))$,
where $C_B(E)$ denotes the centralizer of $E$
in $B$.
Consequently,
\[
    \deg C_B(E)=\frac{\deg B}{[E:k]},
    \quad
    \ind C_B(E)=\ind(B_E).
    \]

\begin{lemma}
\label{embedding}
A finite extension $E/k$ embeds in a central simple
$k$-algebra $B$ if and only if
$[E:k]\cdot\ind(B_E)\mid\deg B$.
\end{lemma}

\begin{proof}
Write $B\simeq M_n(\Delta)$, where $\Delta$
is a central division $k$-algebra.
By \cite[Lem.~2.3]{ShihYangYu}, the field
$E$ embeds in $B$ if and only if
$[E:k]\mid n\cdot c(\Delta_E)$,
where $c(\Delta_E)$ denotes the capacity
(also called the co-index) of $\Delta_E$.
Since
\[
    c(\Delta_E)=\frac{\deg\Delta_E}{\ind(\Delta_E)}
    =\frac{\deg\Delta}{\ind(B_E)},
    \quad
    \deg B=n\deg\Delta,
\]
the assertion follows.
See also \cite[Lem.~7]{Yu2013}.
\end{proof}

\begin{lemma}
\label{field}
Let $L/F$ be a cyclic extension of odd prime degree $p$.
Put $A=M_r(D)$, $A_L=A\otimes_K K_L$, and
$N=\deg A=2^s m r$, and assume $p\nmid N$.
Let $C\in A_L^\times$ have reduced characteristic
polynomial in $K[X]$.
If $R:=K[C]$ is a field, then $d\mid N$ and
$\gcd(d,p)=1$, where $d=[R:K]$.
Moreover,
$\widetilde R:=K_L[C]\simeq R\otimes_K K_L$
is a field, and there exists $P\in A_L^\times$
such that
$
P^{-1}CP\in A.
$
\end{lemma}
\begin{proof}
Write $\mu_C\in K[X]$ for the minimal polynomial of $C$ over $K$, and let $\chi_C\in K[X]$ be the reduced characteristic polynomial of
$C$ in $A_L$; it is monic of degree $N$. Since $K[C]$ is a field, $\mu_C$ is
irreducible over $K$.

We claim that $\chi_C=\mu_C^{a}$ for some $a\ge1$. Fix a splitting
$A_L\otimes_{K_L}\overline{K_L}\simeq M_N(\overline{K_L})$ and let
$\overline C\in M_N(\overline{K_L})$ be the image of $C$. If $\lambda\in\overline{K_L}$ is a root of
$\chi_C$, it is an eigenvalue of $\overline C$; for an eigenvector $v\ne0$,
$\mu_C(\lambda)v=\mu_C(\overline C)v=0$, so $\mu_C(\lambda)=0$.
Now let $f\in K[X]$ be a monic irreducible factor of $\chi_C$ (note
$\chi_C\in K[X]$). A root $\lambda$ of $f$ is a root of $\chi_C$, hence of
$\mu_C$; as $f$ and $\mu_C$ are both monic irreducible over $K$ with the
common root $\lambda$, we get $f=\mu_C$. Hence every irreducible factor of
$\chi_C$ in $K[X]$ equals $\mu_C$, i.e.\ $\chi_C=\mu_C^a$ for some $a\ge1$.
  Taking degrees gives $N=ad$, where $d=\deg\mu_C=[R:K]$. Hence $d\mid N$, and since $p\nmid N$ we get $\gcd(d,p)=1$.

The algebra $R_L=R\otimes_K K_L$ has dimension
$dp$ over $K$.
Let $\mathfrak m$ be a maximal ideal of $R_L$.
Since $R$ and $K_L$ are fields, the natural
maps
$R\longrightarrow R_L/\mathfrak m$ and $K_L\longrightarrow R_L/\mathfrak m$
are injective.
Consequently, both $d$ and $p$ divide
$\dim_K(R_L/\mathfrak m)$.
As $\gcd(d,p)=1$, we obtain
\[
dp\mid\dim_K(R_L/\mathfrak m)
\leq\dim_K R_L=dp.
\]
Thus $\mathfrak m=0$, and hence $R_L$ is a field.
The natural multiplication map
\[
R\otimes_K K_L\longrightarrow A_L,
\quad x\otimes\lambda\longmapsto x\lambda,
\]
is injective, since $R\otimes_K K_L$ is a field.
Its image is precisely $\widetilde R=K_L[C]$.
Set
\[
    i_0=\ind(A\otimes_K R),\quad
    i_1=\ind(A\otimes_K \widetilde R),\quad
    t=N/d.
\]
Both indices divide $\ind(A)=2^s m$.
Since $[\widetilde R:R]=p$, we have
$i_1\mid i_0$ and $i_0\mid p i_1$.
Moreover, $p\nmid \ind(A)$, so $p\nmid i_0$.
Consequently, $i_0=i_1$.

Let $B=C_{A_L}(\widetilde R)$. By the centralizer theorem \cite[Thm. III.5.1(3)]{BO},
we have $A_L\otimes_{K_L}\widetilde R\simeq M_d(B)$.
Hence
\[
\deg B=\frac{N}{d}=t,
\quad
\ind B=\ind(A\otimes_K \widetilde R)=i_1.
\]
Since $\ind B\mid\deg B$ and $i_0=i_1$,
we obtain $i_0\mid t$. Thus
$d\cdot\ind(A\otimes_K R)\mid N$.
By Lemma~\ref{embedding}, there exists
a $K$-embedding $j_0:R\hookrightarrow A$.
Let
\[
j_1:R\otimes_K K_L\xrightarrow{\ \sim\ }\widetilde R\subseteq A_L,\qquad
j_1(x\otimes\lambda)=x\lambda,
\]
be the multiplication isomorphism. Consider the $K_L$-embedding $(j_0)_L:R\otimes_K K_L\longrightarrow A_L$,
defined by
$(j_0)_L(x\otimes\lambda)=j_0(x)\lambda$ for $x\in R$, $\lambda\in K_L$.

Since $R\otimes_K K_L$ is a field and $A_L$ is central
simple over $K_L$, the Skolem--Noether theorem
\cite[Thm.~III.3.2]{BO} yields
$P\in A_L^\times$ such that
\[
P^{-1}j_1(z)P=(j_0)_L(z)
\qquad
(z\in R\otimes_K K_L).
\]
Taking $z=C\otimes 1$, we obtain
$
P^{-1}CP=j_0(C)\in A.
$
\end{proof}

\begin{lemma}
\label{norm}
Let $E/k$ be a cyclic extension of odd prime degree $p$, generated by $\rho$.
Let $H_0$ be central simple over $k$ of degree $t$,
and set $H=H_0\otimes_k E$.
Extend $\rho$ to an automorphism of $H$ via $\id_{H_0}\otimes\rho$.
If $S\in H^\times$ and $c\in k^\times$ satisfy
$
    S\rho(S)\cdots\rho^{p-1}(S)=c
$
and $\gcd(t,p)=1$, then
$
    c\in N_{E/k}(E^\times).
$
\end{lemma}

\begin{proof}
Since $\Nrd_H$ is given by a polynomial with coefficients in $k$ in the
coordinates of $H=H_0\otimes_kE$ with respect to a $k$-basis of $H_0$, and $\rho$
acts on these coordinates while fixing $k$, we have
$\Nrd_H\circ\rho=\rho\circ\Nrd_H$.
Taking reduced norms, we obtain
\[
c^t
=\Nrd_H(c)
 =\prod_{i=0}^{p-1}\Nrd_H(\rho^i(S))=\prod_{i=0}^{p-1}\rho^i(\Nrd_H(S))
 =N_{E/k}(\Nrd_H(S)).
\]
Since $c\in k^\times$, we also have $N_{E/k}(c)=c^p$.
Choose integers $a,b$ such that $at+bp=1$. Then
$N_{E/k}\bigl(\Nrd_H(S)^a c^b\bigr)
=c^{at+bp}=c$.
Thus $c\in N_{E/k}(E^\times)$.
\end{proof}

\begin{lemma}
\label{nor}
Let $E/k$ be a cyclic extension of odd prime degree $p$, generated by
$\rho$.
Let $\theta$ be an automorphism of $E$ of order
dividing $2$, preserving $k$ and commuting with $\rho$.
If
$\theta(c)c=1$ and $c\in N_{E/k}(E^\times)$ for some $c\in k^\times$,
then there is $\zeta\in E^\times$ with
$\theta(\zeta)\zeta=1$ and $N_{E/k}(\zeta)=c$.
\end{lemma}

\begin{proof}
Choose $w\in E^\times$ with $N_{E/k}(w)=c$ and put
$z=w\theta(w)^{-1}$. Since $\theta^2=\id$, we have
$\theta(z)=z^{-1}$.
Moreover, the commutativity of $\theta$ and $\rho$ gives
$N_{E/k}(\theta(w))
=\theta(N_{E/k}(w))
=\theta(c)=c^{-1}.
$
Consequently,
$N_{E/k}(z)=N_{E/k}(w)N_{E/k}(\theta(w))^{-1}=c^2$.
Choose integers $a,b$ such that $2a+pb=1$, and set
$\zeta=z^a c^b$. Using
$N_{E/k}(c)=c^p$, we obtain
$\theta(\zeta)\zeta=1$ and $N_{E/k}(\zeta)=c^{2a+pb}=c$.
\end{proof}

\section{Semilinear isometries and spectral reduction}
\label{sec:spectral}

Fix an admissible cyclic extension $L/F$ of odd prime
degree $p$, and choose a generator $\rho$ of
$\Gal(L/F)$.
We also write $\rho$ for the automorphism
$\id_D\otimes\rho$ of $D_L$.
Its restriction to $K_L$ generates the cyclic
extension $K_L/K$, and $\rho$ commutes with $\tau_L$.

Let $\varphi$ be a nonsingular form over $(D_L,\tau_L)$
with Galois-invariant isometry class, i.e., for every
$\gamma\in\Gal(L/F)$, there exists a
$\gamma$-semilinear isometry of $\varphi$.
Let $V$ be its underlying right $D_L$-module,
and put $r=\dim_{D_L}V$.
Fix a right $D_L$-basis of $V$ and let $\rho$
act coordinatewise on $V$.
Since $\varphi$ is isometric to its $\rho$-twist,
there exists $S\in\GL_r(D_L)$ such that
$T=S\rho$ is a $\rho$-semilinear isometry
of $\varphi$.
Thus
$h(Tx,Ty)=\rho(h(x,y))$
in the hermitian setting, and
$q(Tx)=\rho(q(x))$
in the quadratic setting.
Set $A=M_r(D)$, so that $A_L\simeq M_r(D_L)$,
and
\begin{equation}\label{eq1}
C=T^p=S\rho(S)\cdots\rho^{p-1}(S).
\end{equation}
Define
$\alpha:A_L\longrightarrow A_L$ via $\alpha(X)=S\rho(X)S^{-1}$.

\begin{lemma}
\label{basic}
The operator $C$ is a $D_L$-linear isometry of
$\varphi$, and
$\alpha(C)=C$.
Its minimal polynomial over $K_L$ and its reduced
characteristic polynomial belong to $K[X]$.
The minimal polynomials over $K$ and $K_L$ coincide.
Write $b(x,y)=x^*By$ and
$X^*=\tau_L(X)^t$, and let
$\sigma(X)=B^{-1}X^*B$
be the adjoint involution with respect to $b$.
Then
$\sigma(C)=C^{-1}$ and
$\alpha\sigma=\sigma\alpha$.
\end{lemma}

\begin{proof}
Since $T$ is a semilinear isometry, iterating $p$ times
shows that $C=T^p$ is a linear isometry of $\varphi$.
In particular, $C$ preserves the associated hermitian
form $b$, and in the quadratic setting it also preserves $q$.

Applying $\rho$ to \eqref{eq1} and then multiplying by $S$ on the left gives $S\rho(C)=CS$.
Hence $\alpha(C)=S\rho(C)S^{-1}=C$. Therefore $C$ and $\rho(C)$ are conjugate,
so their minimal and reduced characteristic polynomials are fixed
by $\rho$. Since $\rho$ generates $\Gal(K_L/K)$, these
polynomials lie in $K[X]$. In particular,
$\mu_{C,K_L}\in K[X]$, hence it annihilates $C$ over $K$;
together with $\mu_{C,K_L}\mid\mu_{C,K}$ this gives
$\mu_{C,K}=\mu_{C,K_L}$.

Write
$b(x,y)=x^*By$. Note that $b(Tx,Ty)=\rho\bigl(b(x,y)\bigr)$ for all $x,y\in V$.
Using $T=S\rho$ and $b(x,y)=x^*By$, this becomes
\[
\rho(x)^*S^*BS\rho(y)
=\rho(x)^*\rho(B)\rho(y).
\]
Since $b$ is nonsingular, it follows that
$S^*BS=\rho(B)$.
Moreover,
\[
x^*C^*BCy=b(Cx,Cy)=b(x,y)=x^*By,
\]
so $
C^*BC=B.
$
Hence,
$\sigma(C)=B^{-1}C^*B=C^{-1}$.
Finally, using $\rho(X^*)=\rho(X)^*$ and $S^*BS=\rho(B)$, we obtain
\[\alpha(\sigma(X))
    =S\rho(B)^{-1}\rho(X)^*\rho(B)S^{-1}
    =B^{-1}(S^*)^{-1}\rho(X)^*S^*B
    =\sigma(\alpha(X)).\qedhere\]
\end{proof}

\begin{lemma}
\label{her}
If $\varphi=h$ is anisotropic hermitian, then
$K[C]$ is reduced.
\end{lemma}

\begin{proof}
The commutative algebra $R=K[C]$ is preserved by
$\sigma$, since $\sigma(C)=C^{-1}\in R$ by Lemma \ref{basic}.
Its nilradical $J$ is also preserved.
If $J\ne0$, choose $\nu\geq1$ maximal with
$J^\nu\ne0$ and choose $0\ne X\in J^\nu$.
Then
$\sigma(X)X\in J^{2\nu}=0$.
Choose $v\in V$ such that $Xv\ne0$. Then
$h(Xv,Xv)=h(v,\sigma(X)Xv)=0$,
contradicting anisotropy.
\end{proof}

\begin{lemma}\label{frobenius}
In the quadratic setting there is $e\ge0$ such that, after replacing
$T,\rho,C$ by $T^{2^e},\rho^{2^e},C^{2^e}$, the following hold:
$\rho$ still generates $\Gal(L/F)$; $T=S\rho$ with $S\in\GL_r(D_L)$ is
a $\rho$-semilinear isometry of $q$ with $T^p=C$; and the minimal polynomial of $C$ over $F$ is separable; in particular, $F[C]$ is
reduced.
\end{lemma}

\begin{proof}
Recall that $\car F=2$, $K=F$ and $N=\deg A=r2^s$.
Let $\bar L$ be an algebraic closure of $L$ and fix an isomorphism
$A_L\otimes_L\bar L\simeq M_N(\bar L)$; regard $C$ as an element of
$M_N(\bar L)$. Since $A_L$ is finite-dimensional over $F$, $F[C]$ is
finite-dimensional, so every eigenvalue of $C$ is a root of the minimal
polynomial of $C$ over $F$ and hence algebraic over $F$.

Let $\lambda_1,\dots,\lambda_k$ be the distinct eigenvalues of $C$ and
$V_{\bar L}=\bigoplus_i V_i$ the decomposition into generalized eigenspaces.
On $V_i$ write $C|_{V_i}=\lambda_iI+Z_i$ with $Z_i$ nilpotent; since
$\dim V_i\le N$, $Z_i^N=0$. Let $f_i\in F[X]$ be the minimal polynomial of
$\lambda_i$ over $F$; as $\car F=2$,
$f_i(X)=g_i(X^{2^{n_i}})$ with $g_i\in F[X]$ irreducible separable.
Choose $e$ with $2^e\ge N$ and $e\ge n_i$ for all $i$.
Then $\lambda_i^{2^{n_i}}$ is a root of $g_i$, hence separable over $F$, and
so is $\lambda_i^{2^e}$. Since $\car F=2$ and
$\lambda_iI$ commutes with $Z_i$,
\[
C^{2^e}|_{V_i}=(\lambda_iI+Z_i)^{2^e}=\lambda_i^{2^e}I+Z_i^{2^e}=\lambda_i^{2^e}I.
\]
Thus $C'=C^{2^e}$ is diagonalizable over $\bar L$ with eigenvalues
$\lambda_i^{2^e}$, each separable over $F$. Let $g$ be the least common
multiple of the minimal polynomials over $F$ of the $\lambda_i^{2^e}$; it
is separable, and $g(C')$ is diagonalizable with all eigenvalues $0$, so
$g(C')=0$. Thus $\mu_{C',F}$ is separable, and hence $F[C']\simeq F[X]/(\mu_{C',F})$ is a
finite product of fields; in particular, it is reduced.

Since $p$ is odd, $\gcd(2^e,p)=1$ and $\rho'=\rho^{2^e}$ generates $\Gal(L/F)$.
Put $T'=T^{2^e}$. From $q(Tx)=\rho(q(x))$ (where $\rho$ acts on
$\overline{D_L}=D_L/\Alt(D_L,\tau_L)$) we get
$q(T'x)=\rho'(q(x))$, and $T'=S'\rho'$ with
$S'=S\rho(S)\cdots\rho^{2^e-1}(S)\in\GL_r(D_L)$. The polar form is
determined by $q$ (Lemma~\ref{polar}), so $T'$ is a semilinear isometry of
$q$. Finally $(T')^p=(T^p)^{2^e}=C'$.
\end{proof}

\begin{lemma}\label{fs}
Let $q$ be a generalized quadratic form in the
quadratic setting, and suppose $F[C]$ is reduced.
There is a $T$-stable orthogonal decomposition
\[
    V=V_0\perp V_*,
    \qquad
    V_0=\ker(C-1),\quad
    V_*=(C-1)V.
\]
Here either summand is allowed to be zero.
The restriction of $C$ to $V_0$ is the identity,
and $C-1$ is invertible on $V_*$.
Moreover, the restricted polar forms on both
summands are nonsingular.
\end{lemma}

\begin{proof}
Since $F[C]\simeq F[X]/(\mu_{C,F})$ is reduced,
the minimal polynomial $\mu_{C,F}$ has no
repeated irreducible factors.
If $X-1$ does not occur in $\mu_{C,F}$, then $C-1$ is invertible and the result is evident.
Otherwise, it occurs only
once and is coprime to the remaining factors. By the Chinese remainder
theorem,
\[
    V=\ker(C-1)\oplus(C-1)V,
\]
and $C-1$ is invertible on $(C-1)V$.
Since $C$ preserves $b$, for all $x,y\in V$ we have
$b(C^{-1}x,y)=b(x,Cy)$.
For $x\in V_0=\ker(C-1)$ and $y\in V$, we have
\[
b(x,(C-1)y)=b(x,Cy)-b(x,y)=b(C^{-1}x,y)-b(x,y)={b((C^{-1}-1)x,y)}=0.
\]
Thus $V_0$ and $V_*=(C-1)V$ are orthogonal for $b$, hence also
for $q$.

Since $V=V_0\oplus V_*$ is orthogonal and $b$ is nonsingular,
$b|_{V_0},\ b|_{V_*}$ are nonsingular.
Finally, $T$ commutes with $C=T^p$, so it preserves both
$\ker(C-1)$ and $(C-1)V$.
\end{proof}

\begin{lemma}
\label{primitive}
Suppose $R=K[C]$ is reduced.
In the hermitian setting, assume the form is
anisotropic.
In the quadratic setting, assume the form is
anisotropic and $C-1$ is invertible.
Then every primitive idempotent of $R$ is fixed by
$\sigma$ and $\alpha$.
Consequently the primitive idempotents give a
$T$-stable orthogonal decomposition into anisotropic
nonsingular forms, on each of which the algebra
generated by the restricted $C$ is a field.
\end{lemma}

\begin{proof}
By Lemma \ref{basic}, we have $\sigma(C)=C^{-1}$ and $\alpha(C)=C$.
Hence,
$\sigma$ preserves $R$ and $\alpha$ fixes $R$
pointwise.
Let $e$ be a primitive idempotent.
since $\sigma$ preserves $R$, $\sigma(e)$ is also a primitive idempotent.
Suppose $\sigma(e)\ne e$.
Since $R$ is reduced and commutative, distinct primitive idempotents
are orthogonal, so $\sigma(e)e=0$. It follows that
$b(ex,ey)=b(x,\sigma(e)ey)=0$
for all $x,y\in V$, so $b$ vanishes identically on
$eV$.
In the hermitian setting, this contradicts anisotropy.
In the quadratic setting, this implies that $q|_{eV}$ is
additive. Since $C$ preserves $q$,
\begin{equation}\label{eq2}
q((C-1)x)=q(Cx)+q(x)+\overline{b(Cx,x)}=q(x)+q(x)+\overline{b(Cx,x)}=\overline{b(Cx,x)},
\end{equation}
for all $x\in eV$.
Since $e\in R=K[C]$, it commutes with $C$.
Thus $Cx=Cex=eCx\in eV$ for every $x\in eV$.
It follows from \eqref{eq2} that $q((C-1)x)=0$ for all $x\in eV$.
But $C-1$ is invertible on $V$ and $eV$ is
$C$-stable, so $C-1$ restricts to an automorphism of $eV$.
Thus $q$ vanishes on all of $eV$, contradicting anisotropy. Therefore $\sigma(e)=e$, so the summands $eV$ and $(1-e)V$ are orthogonal
for $b$. The restricted forms are nonsingular and remain anisotropic.
Also, since $\alpha(e)=e$, we have $S\rho(e)S^{-1}=e$, hence $S\rho(e)=eS$, which implies that
\[T(ex)=S\rho(ex)=S\rho(e)\rho(x)=eS\rho(x)=eT(x),\]
 for all $x\in V$.
Hence $T$ preserves each summand. Finally, the restriction map $eR\rightarrow\End_{D_L}(eV)$
is injective, since $eR$ is a field and its identity $e$
acts as the identity on the nonzero module $eV$.
Its image is precisely the $K$-algebra generated by
$C|_{eV}$.
\end{proof}

\begin{corollary}
\label{spectral}
Let $\varphi$ be anisotropic with Galois-invariant
isometry class.
In the hermitian setting, $\varphi$ is a $T$-stable
orthogonal sum of anisotropic nonsingular forms
whose restricted $K[C]$ are fields.
In the quadratic setting, after the replacement in
Lemma~\ref{frobenius}, $\varphi$ is a $T$-stable
orthogonal sum consisting of:
\begin{enumerate}
\item a possibly zero anisotropic nonsingular summand
      on which $C=1$;
\item anisotropic nonsingular summands on which
      $F[C]$ is a field and $C-1$ is invertible.
\end{enumerate}
\end{corollary}

\begin{proof}
For hermitian forms, Lemma~\ref{her} gives that
$K[C]$ is reduced. Since $b$ is nonsingular and $h$ is anisotropic,
Lemma~\ref{primitive} applies to $V$ and yields the desired
$T$-stable orthogonal decomposition into anisotropic nonsingular
summands on which $K[C]$ is a field.

For quadratic forms, first apply Lemmas~\ref{frobenius} and
\ref{fs}, after replacing $C$ by $C^{2^e}$. We obtain
\[
    V=V_0\perp V_*,
    \quad
    V_0=\ker(C-1),\quad
    V_*=(C-1)V,
\]
where $C=1$ on $V_0$.
Since $b$ is nonsingular and the decomposition is
orthogonal, the restricted polar forms on both
$V_0$ and $V_*$ are nonsingular.
 Both summands are $T$-stable because $T$ commutes with
$C$.
Now apply Lemma~\ref{primitive} to $V_*$. It gives a
$T$-stable orthogonal decomposition of $V_*$ into anisotropic
nonsingular summands on which $F[C]$ is a field and $C-1$ is
invertible. Together with $V_0$, this gives the required
decomposition.
\end{proof}

\section{Descent of invariant anisotropic forms}
\label{sec:descent}

We continue with the notation of the preceding sections.
In particular, $\varphi$ denotes a nonsingular
form of either type over $(D_L,\tau_L)$.
We show that if $\varphi$ is anisotropic and its
isometry class is Galois-invariant, then it descends
to $(D,\tau)$ whenever its rank is sufficiently small
relative to $[L:F]$.

\begin{lemma}
\label{descent}
Suppose $U$ is a
$D_L$-linear isometry of $\varphi$ satisfying
\[
    U\alpha(U)\cdots\alpha^{p-1}(U)=C^{-1}.
\]
Then $\varphi$ descends to a nonsingular form over $(D,\tau)$ of
the same type and rank.
The associated hermitian form descends as well.
If $\varphi$ is anisotropic, so is its descended form.
\end{lemma}

\begin{proof}
Set $T'=UT$. Since $U$ is $D_L$-linear and $T$ is
$\rho$-semilinear, $T'$ is $\rho$-semilinear. Both $U$ and $T$
are isometries of $\varphi$, hence so is $T'$.
We claim that $(T')^p=1$. Since $T=S\rho$, and $\rho X\rho^{-1}$ is the entrywise
image of $X$ under $\rho$, we have
$TXT^{-1}=S\rho(X)S^{-1}=\alpha(X)$
for all $X\in\End_{D_L}(V)$. Induction gives
\[
    (UT)^p=U\,\alpha(U)\cdots\alpha^{p-1}(U)\,T^p.
\]
Hence $(T')^p=1$ and the rule $\rho^i\cdot x=(T')^i x$ for $0\leq i<p$
defines a semilinear action of $\Gal(L/F)$
on $V$. By Galois descent,
$V_0=V^{T'}$ is a right $D$-vector space and the natural map
$V_0\otimes_F L\to V$ is an isomorphism of right
$D_L$-modules.

For $x,y\in V_0$, semilinear invariance of the associated
hermitian form $b=b_\varphi$ implies
\[
\rho(b(x,y))=b(T'x,T'y)=b(x,y).
\]
Thus $b(x,y)\in(D_L)^{\Gal(L/F)}=D$, and
$b$ restricts to a hermitian form $b_0$ over $(D,\tau)$.

In the quadratic setting, scalar extension gives an isomorphism
of $L$-vector spaces
\[
\overline D_L
= D_L/\Alt(D_L,\tau_L)
\simeq (D/\Alt(D,\tau))\otimes_F L
=\overline D\otimes_F L.
\]
The isomorphism is $\Gal(L/F)$-equivariant; hence
\[
(\overline D_L)^{\Gal(L/F)}
\simeq (\overline D\otimes_F L)^{\Gal(L/F)}
\simeq\overline D.
\]
For every $x\in V_0$, semilinear invariance gives
$\rho(q(x))=q(T'x)=q(x)$.
Consequently $q(x)$ belongs to the image of $\overline D$
inside $\overline D_L$. It therefore defines a map
$q_0:V_0\to\overline D$, whose polar form is $b_0$.
The defining identities restrict to $V_0$, and extension of
scalars recovers $\varphi$.

Nonsingularity descends by faithful
scalar extension. Finally, an isotropic vector of the
descended form would remain isotropic over $L$, contradicting
the anisotropy of $\varphi$.
\end{proof}

\begin{corollary}
\label{fd}
A form with a semilinear isometry $T$ satisfying
$T^p=1$ descends.
\end{corollary}

\begin{proof}
Use Lemma~\ref{descent} with $U=1$ and $C=1$.
\end{proof}

\begin{lemma}
\label{qu}
Assume $\car F=2$. Let $q_1,q_2$ be generalized quadratic forms
on $V$ with the same polar form. Suppose a $D$-linear automorphism
$c$ of $V$ preserves both $q_1$ and $q_2$, and that $c-1$ is
invertible. Then $q_1=q_2$.
Consequently, for the fixed form $q$, if $c$ preserves $q$,
$c-1$ is invertible, and $u$ is a polar isometry of $q$ commuting
with $c$, then $u$ preserves $q$, i.e. $u$ is an isometry of $q$.
\end{lemma}

\begin{proof}
The difference $\delta=q_1-q_2$ has zero polar form, hence
$\delta(x+y)=\delta(x)+\delta(y)$
for all $x,y\in V$.
Since $c$ preserves both $q_1$ and $q_2$, we have
$\delta(cx)=\delta(x)$. Therefore
\[
    \delta((c-1)x)=\delta(cx+x)=\delta(cx)+\delta(x)=0.
\]
Since $c-1$ is invertible, every $y\in V$ is of the form
$y=(c-1)x$, so  $\delta=0$, i.e. $q_1=q_2$.

For the final assertion, apply the first part to $q$ and
$q_u(x)=q(ux)$. They have the same polar form because $u$ is a
polar isometry. Moreover, $c$ preserves both:
\[
    q_u(cx)=q(ucx)=q(cux)=q(ux)=q_u(x),
\]
since $u$ commutes with $c$. Hence $q_u=q$, so $u$ preserves
$q$.
\end{proof}

\begin{proposition}
\label{fco}
Let $\varphi$ be a form of either type over $(D_L,\tau_L)$,
of positive rank $r<p$, and with
Galois-invariant isometry class.
Write $V$ for the underlying $D_L$-module,
put $A=M_r(D)$.
In the quadratic setting, use the modified
semilinear data $(T,\rho,C)$ from
Lemma~\ref{frobenius}, with $S$ and
$\alpha$ defined accordingly.
Suppose $K[C]$ is a field.
Then $\varphi$ descends to $(D,\tau)$.
\end{proposition}

\begin{proof}
Since $K[C]$ is a field, either $C=1$
or $C-1$ is invertible.
If $C=1$, then $T^p=1$, so
Corollary~\ref{fd} applies.
Assume henceforth $C\ne1$.

Put $N=\deg A=2^s m r$, $d=[K[C]:K]$, and $t=N/d$.
Since $L/F$ is admissible, $p\nmid m$.
Moreover, $r<p$ gives $p\nmid r$.
Consequently, $p\nmid N$.
By Lemmas~\ref{basic} and~\ref{field}, after a suitable
coordinate change (replace $S$, $B$, $C$ by
$S'=P^{-1}S\rho(P)$, $B'=P^*BP$, $C'=P^{-1}CP$ with
$P\in A_L^\times$), the element $C'$ lies in $A$.
Define the transformed maps
\[
\alpha'(X)=S'\rho(X)S'^{-1},\qquad
\sigma'(X)=(B')^{-1}X^*(B').
\]
Then
\[
\alpha'(P^{-1}XP)=P^{-1}\alpha(X)P,\qquad
\sigma'(P^{-1}XP)=P^{-1}\sigma(X)P,
\]
where $\sigma$ is the adjoint involution with respect to $b$.
In particular, $\alpha'$ and $\sigma'$ commute.
We define the transformed hermitian form by
$b'(x,y):=b(Px,Py)$,
so that its matrix with respect to the chosen basis is $B'=P^*BP$.
In the hermitian setting we put $\varphi':=b'$.
In the quadratic setting, let $q'$ be the transformed form
$q'(x):=q(Px)$,
and put $\varphi':=(q',b')$, where $b'=b_{q'}$.
Moreover, $\varphi'$ is isometric to $\varphi$, and $T':=P^{-1}TP=S'\rho$
is a $\rho$-semilinear isometry of $\varphi'$.
Set
\[
R=K[C'],\quad \widetilde R=K_L[C'],\quad H_0=C_A(R),\quad H=C_{A_L}(\widetilde R).
\]
Lemma~\ref{field} shows that $d\mid N$,
$\widetilde R\simeq R\otimes_K K_L$ is a field, and $C'\in R^\times$.
Since $p\nmid d$, the fields $R$ and $K_L$ are linearly disjoint
over $K$; hence $\widetilde R/R$ is cyclic of degree $p$, generated by
$\rho_{\widetilde R}=\id_R\otimes\rho$, and $\gcd(t,p)=1$.

By the centralizer theorem \cite[Thm. III.5.1(1)]{BO},
$H_0$ is central simple over $R$ of degree $t$,
and $H$ is central simple over $\widetilde R$ of the
same degree.
Since centralizers commute with flat scalar
extension, we have
$H=H_0\otimes_K K_L=H_0\otimes_R \widetilde R$.

Since $C'\in A$ and $\rho$ acts trivially on $A$, we have $\rho(C')=C'$.
Since $\alpha(C)=C$ by Lemma \ref{basic}, the transformation formula for $\alpha'$ gives
\[
\alpha'(C')=\alpha'(P^{-1}CP)=P^{-1}\alpha(C)P=C'.
\]
Hence
$
S'C'S'^{-1}=\alpha'(C')=C',
$
so $S'C'=C'S'$. As $K_L$ is central in $A_L$, it follows that $S'$
centralizes $\widetilde R=K_L[C']$. Thus $S'\in H^\times$. The coefficientwise
action of $\rho$ on $H$ is $\id_{H_0}\otimes\rho_{\widetilde R}$, under the
tensor-product identification above.

The equality ($T')^p=P^{-1}T^pP=P^{-1}CP=C'$ gives
\[
S'\rho(S')\cdots\rho^{p-1}(S')=C'.
\]
Since $\widetilde R/R$ is cyclic, Lemma~\ref{norm} yields $C'\in N_{\widetilde R/R}(\widetilde R^\times)$.

Since $\sigma'$ restricts to $\tau_L$ on $K_L$ and
$\sigma'(C')=C^{\prime -1}$, the involution $\sigma'$
preserves $\widetilde R=K_L[C']$ and $R=K[C']$.
The restriction $\theta=\sigma'|_{\widetilde R}$ is an automorphism of $\widetilde R$ of order
dividing $2$, preserves $R$, satisfies $\theta(C')=C^{\prime -1}$, and hence
$\theta(C')C'=1$.
Since $S'$ centralizes $\widetilde R$, we have $\alpha'|_{\widetilde R}=\rho_{\widetilde R}$; the commutation
$\alpha'\sigma'=\sigma'\alpha'$ then implies that $\theta$ commutes with
$\rho_{\widetilde R}$.
Applying Lemma~\ref{nor}, we obtain
$\zeta\in \widetilde R^\times$ with
$\sigma'(\zeta)\zeta=1$
and
$N_{\widetilde R/R}(\zeta)=C'$.
In particular, $\zeta$ preserves $b'$, because
\[
    b'(\zeta x,\zeta y)
    =b'(x,\sigma'(\zeta)\zeta y)
    =b'(x,y).
\]
In the hermitian setting, this directly says that $\zeta$ is an
isometry of $\varphi'$.

In the quadratic setting, $C'$ preserves $q'$,
since $C'=(T')^p$ and $T'$ is a $\rho$-semilinear isometry of $q'$.
Moreover, $\zeta$ commutes with $C'$, as both belong to
$\widetilde R$. Since $C'-1$ is invertible, Lemma~\ref{qu}
shows that $\zeta$ preserves $q'$ itself.

Set $U=\zeta^{-1}$. Then $U$ is an isometry of $\varphi'$, and
\[
U\alpha'(U)\cdots\alpha'^{p-1}(U)
=N_{\widetilde R/R}(\zeta)^{-1}=C^{\prime -1}.
\]
Lemma~\ref{descent} now applies to the transformed data
$(\alpha',U,C')$: it descends the transformed form $\varphi'$. Since $\varphi'$ is isometric to $\varphi$ (as noted above), the descended form over $(D,\tau)$ has scalar extension to $L$ that is isometric to $\varphi$; hence $\varphi$ descends to $(D,\tau)$.
\end{proof}

\begin{theorem}
\label{invariant}
Let $L/F$ be an admissible cyclic extension of odd
prime degree $p$.
Every anisotropic form $\varphi$ of either type over
$(D_L,\tau_L)$, of rank less than $p$ and with
Galois-invariant isometry class, descends to an
anisotropic form over $(D,\tau)$.
\end{theorem}

\begin{proof}
The zero form descends trivially, so assume $\rank\varphi>0$.

In the quadratic setting, apply the replacement of
Lemma~\ref{frobenius}: replace $T$, $\rho$, $C$ by
$T^{2^e}$, $\rho^{2^e}$, $C^{2^e}$ for a sufficiently large $e$.
We keep the notation $T,\rho,C$ for the replaced
objects.
By Corollary~\ref{spectral}, $\varphi$ admits a $T$-stable
orthogonal decomposition into summands of the following two types.
\\
(1) A possibly zero summand on which $C=1$. On this summand
$T^p=C=1$, so Corollary~\ref{fd} applies.
\\
(2) Anisotropic nonsingular summands on which the algebra
generated by the restricted $C$ is a field, and, in the quadratic
setting, $C-1$ is invertible. Each such summand has positive rank,
and its isometry class is Galois-invariant because
$T$ preserves it and is a semilinear isometry.
More explicitly, writing
$r_\varphi=\rank \varphi<p$,
let $V_i$ be any nonzero summand of this second type, and put
$
r_i=\dim_{D_L}V_i.
$
Since the decomposition is a direct sum, we have
$
1\leq r_i\leq r_\varphi<p.
$
Consequently, $p\nmid r_i$. Since $L/F$ is admissible,
$p\nmid m$, and since $p$ is odd, $p\nmid 2^s$.
Therefore,
\[
p\nmid N_i=2^s m r_i=\deg M_{r_i}(D).
\]
Thus, the hypothesis of Lemma~\ref{field} is satisfied
for each such summand in the proof of
Proposition~\ref{fco}.
Hence Proposition~\ref{fco} applies
and descends each of these summands.

Taking the orthogonal sum of the descended forms gives
a form $\psi$ over $(D,\tau)$ whose scalar extension
to $L$ is isometric to $\varphi$.
Since scalar extension is faithful, $\psi$ is
anisotropic.
\end{proof}

\section{The main result}
\label{sec:support}

We continue with the notation of the preceding sections.
In this section, however, $\varphi$ is a nonsingular
form defined over $(D,\tau)$, rather than over
$(D_L,\tau_L)$ as in the preceding sections.
Its scalar extension to $(D_L,\tau_L)$ is denoted
by $\varphi_L$.
We prove the main theorem by combining the descent
result of Theorem~\ref{invariant} with a bounded-rank
reduction.

\begin{proposition}
\label{bounded}
Let $L/F$ be an admissible cyclic extension of odd
prime degree $p$, and let $\varphi$ be a form of either
type over $(D,\tau)$.
If $\varphi$ is anisotropic and
$\rank\varphi\leq p+1$,
then $\varphi_L$ is anisotropic.
\end{proposition}

\begin{proof}
Suppose $\varphi_L$ is isotropic.
Take a Witt decomposition
\[
    \varphi_L\simeq \varphi_{\mathrm{met}}\perp\varphi_{\mathrm{an}},
    \quad
    \rank \varphi_{\mathrm{met}}=2j,\quad j\geq1,
\]
where $\varphi_{\mathrm{met}}$ is metabolic,
$\varphi_{\mathrm{an}}$ is anisotropic or zero with
$\rank\varphi_{\mathrm{an}}<p$.
Since $\varphi_L$ is obtained by scalar extension
from $F$, its Witt class is Galois invariant.
Galois twisting preserves Witt decompositions,
so $\gamma(\varphi_{\mathrm{an}})$ and
$\varphi_{\mathrm{an}}$ are anisotropic representatives
of the same Witt class (see Proposition~\ref{witt}).
Theorem~\ref{invariant} gives an anisotropic
descended form $\varphi_0$, including the zero case.
It is nonsingular because its scalar extension is
nonsingular.
Now $\varphi_L$ and $(\varphi_0)_L$ are Witt-equivalent, so the restriction map $\res_{L/F}$ sends the Witt classes of $\varphi$ and $\varphi_0$ to the same element. Proposition~\ref{injectivity} then shows that $\varphi$ and $\varphi_0$ are Witt-equivalent over $F$. Uniqueness of anisotropic representatives then gives $\varphi\simeq\varphi_0$, contradicting $\rank\varphi_0=\rank\varphi-2j$.
\end{proof}

\begin{lemma}
\label{support}
Let $E/F$ be a finite extension of odd degree $n$,
and let $\varphi$ be a form of either type over $(D,\tau)$.
If $\varphi$ is anisotropic and $\varphi_E$ is
isotropic, there is a right $D$-subspace $W$ such that
$1\leq\dim_D W\leq
\min\{n,\rank\varphi\}$,
the restriction $\varphi|_W$ is anisotropic, and
$(\varphi|_W)_E$ is isotropic.
In the hermitian setting the restriction is
automatically nonsingular.
\end{lemma}

\begin{proof}
Choose an $F$-basis $a_1,\ldots,a_n$ of $E$ and write
a nonzero isotropic vector as
$v=\sum_{i=1}^n v_i\otimes a_i$.
Let $W$ be the right $D$-span of the $v_i$.
It has the required dimension, its restriction is
anisotropic, and $v\in W_E$ witnesses isotropy over
$E$.
For a hermitian restriction, a nonzero vector in its
radical would be isotropic, so its radical is zero.
\end{proof}

\begin{theorem}
\label{cyclic}
Let $L/F$ be an admissible cyclic extension
of odd prime degree.
Every anisotropic form of either type over
$(D,\tau)$, with nonsingular polar form in the quadratic setting, remains anisotropic over $(D_L,\tau_L)$.

\end{theorem}

\begin{proof}
Write $p=[L:F]$.
Suppose $\varphi$ is anisotropic but $\varphi_L$ is
isotropic.
By Lemma~\ref{support}, there is an anisotropic
restriction $\varphi'$ of rank
$1\leq r\leq p$
whose extension to $L$ is isotropic.

In the hermitian setting, $\varphi'$ is nonsingular.
Proposition~\ref{bounded} immediately gives a
contradiction.
In the quadratic setting with $D\ne F$, apply
Proposition~\ref{reg} to $\varphi'$.
It gives a nonsingular anisotropic form
$\varphi'_t$ of rank $r$ over $D((t))$ which becomes
isotropic over $L((t))$.
But
$L((t))/F((t))$
is cyclic of degree $p$, and $D((t))$ is central
division of degree $2^s$.
Since $r\leq p<p+1$, Proposition~\ref{bounded},
applied over $F((t))$, gives a contradiction.

Finally, in the quadratic setting with $D=F$,
generalized quadratic forms are ordinary quadratic
forms in characteristic $2$.
Preservation of anisotropy under odd-degree
extensions, including singular forms is proved in
\cite[Cor.~18.5]{EKM}.
\end{proof}

\begin{theorem}
\label{galois}
Let $E/F$ be an admissible finite Galois extension.
Every anisotropic form $\varphi$ of either type, with nonsingular polar form in the quadratic setting, remains
anisotropic over $E$.
Consequently,
$i(\varphi_E)=i(\varphi)$ for every nonsingular
form $\varphi$.
\end{theorem}

\begin{proof}
By the Feit--Thompson theorem \cite{FT}, the group
$G=\Gal(E/F)$ is solvable.
Choose a composition series
\[
    G=G_0\triangleright G_1\triangleright\cdots
      \triangleright G_\ell=\{1\}
\]
with cyclic quotients of odd prime order.
The fixed fields give a tower
\[
    F=F_0\subset F_1\subset\cdots\subset F_\ell=E
\]
whose successive extensions are cyclic of odd prime
degree.
For each $i$, the degree $[F_i:F]$ divides
$[E:F]$ and is therefore coprime to $m$.
The index divisibilities imply that
$D_{F_i}$ remains a division algebra of degree
$2^s m$. Hence every successive extension
$F_{i+1}/F_i$ is admissible for
$(D_{F_i},\tau_{F_i})$.
Thus Theorem~\ref{cyclic} applies at
every step of the tower.
By induction on $i$, $\varphi_{F_i}$ is anisotropic for all $i$.

Moreover, let $\varphi$ be nonsingular. By
Proposition~\ref{witt}, write
$\varphi\simeq\varphi_{\mathrm{met}}\perp\varphi_{\mathrm{an}}$
 and
$\rank \varphi_{\mathrm{met}}=2i(\varphi)$,
with $\varphi_{\mathrm{met}}$ metabolic and $\varphi_{\mathrm{an}}$
anisotropic. Since metabolicity is preserved under scalar
extension and $(\varphi_{\mathrm{an}})_E$ remains
anisotropic, we obtain
$
i(\varphi_E)=i((\varphi_{\mathrm{met}})_E)
=\frac12\rank \varphi_{\mathrm{met}}
=i(\varphi).
$
\end{proof}

\begin{corollary}
\label{csa-galois}
Problem~\ref{csa} has an
affirmative answer for a finite Galois extension $E/F$ of odd degree with
$\gcd([E:F],\ind A)=1$: every anisotropic form of either type
over $(A,\sigma)$ remains anisotropic after scalar
extension to $(A_E,\sigma_E)$.
\end{corollary}

\begin{proof}
By Theorem~\ref{mor} and Remark \ref{new}, it suffices to establish the division
case. This is precisely Theorem~\ref{galois}, applied to the
division algebra $D$ Brauer equivalent to $A$ and to the
involution $\tau$ provided by Lemma~\ref{md}. The index
condition is preserved because $\ind D=\ind A$.
\end{proof}

\section*{Statements and Declarations}

\noindent\textbf{Funding.}
No funding was received for conducting this study.

\medskip

\noindent\textbf{Competing Interests.}
The author declares that there are no competing interests.

\medskip

\noindent\textbf{Use of AI-assisted Tools.}
During the preparation of this manuscript, ChatGPT
(OpenAI) was used to assist with drafting and
restructuring the abstract, language editing,
LaTeX formatting, and preliminary bibliographic
checks. The author reviewed the suggested changes
and takes full responsibility for the accuracy,
originality, and integrity of the manuscript.

\end{document}